\documentclass[11pt]{amsart}

\usepackage[english]{babel}
\usepackage[T1]{fontenc}
\usepackage[utf8]{inputenc}

\usepackage{mathpazo}
\usepackage{graphicx}
\usepackage{microtype}
\usepackage{enumerate}
\usepackage{csquotes}
\usepackage[dvipsnames]{xcolor}
\usepackage{url}

\usepackage{breakurl}
\usepackage[breaklinks,bookmarks, colorlinks, citecolor=MidnightBlue!75!black, urlcolor=MidnightBlue!75!black, linkcolor=MidnightBlue!75!black]{hyperref}

\calclayout

\AtBeginDocument{%
   \def\MR#1{}
}

\usepackage{amsmath, amssymb, amsthm, amscd}
\usepackage{commath, mathtools, mathrsfs}
\usepackage{thmtools}
\usepackage{tikz, tikz-cd}
\usetikzlibrary{matrix, arrows, decorations.pathmorphing}

\theoremstyle{plain}

\newtheorem{theorem}{Theorem}
\newtheorem{lemma}[theorem]{Lemma}

\newtheorem*{theorem*}{Theorem}
\newtheorem*{lemma*}{Lemma}
\newtheorem*{proposition*}{Proposition}
\newtheorem*{corollary*}{Corollary}

\declaretheorem[numbered=no,name=Theorem A]{TA}
\declaretheorem[numbered=no,name=Theorem B]{TB}
\declaretheorem[numbered=no,name=Theorem C]{TC}

\theoremstyle{definition}

\newtheorem*{conjecture*}{Conjecture}
\newtheorem*{remark*}{Remark}
\newtheorem*{definition*}{Definition}
\newtheorem*{observation*}{Observation}

\newcommand{\on}{\operatorname}

\newcommand{\A}{\mathbb{A}}

\newcommand{\gal}{\on{Gal}}

\newcommand{\Q}{\mathbb{Q}}

\newcommand{\s}{\subseteq}

\newcommand{\spec}{\on{Spec}}

\newcommand{\xar}{\xrightarrow}

\renewcommand{\epsilon}{\varepsilon}

\renewcommand{\P}{\mathbb{P}}
\renewcommand{\phi}{\varphi}

\renewcommand{\tilde}{\widetilde}

\usepackage[giveninits=true,style=alphabetic,backend=biber,maxnames=99,maxalphanames=5,giveninits=true]{biblatex}
\bibliography{blfg}

\author{Giulio Bresciani}
\address[G. Bresciani]{CRM Ennio de Giorgi, Collegio Puteano, Office 21, Piazza dei Cavalieri 3, 56126 Pisa}
\email{giulio.bresciani@sns.it}
\thanks{The author was partially supported by the DFG Priority Program "Homotopy Theory and Algebraic Geometry" SPP 1786}

\title{On the Bombieri-Lang Conjecture over finitely generated fields}
\date{}

\begin{document}

\begin{abstract}
	The \emph{strong Bombieri-Lang conjecture} postulates that, for every variety $X$ of general type over a field $k$ finitely generated over $\mathbb{Q}$, there exists an open subset $U\subset X$ such that $U(K)$ is finite for every finitely generated extension $K/k$. The \emph{weak Bombieri-Lang conjecture} postulates that, for every positive dimensional variety $X$ of general type over a field $k$ finitely generated over $\mathbb{Q}$, the rational points $X(k)$ are not dense. Furthermore, Lang conjectured that every variety of general type $X$ over a field of characteristic $0$ contains an open subset $U\subset X$ such that every subvariety of $U$ is of general type, this statement is usually called \emph{geometric Lang conjecture}.
	
	We reduce the strong Bombieri-Lang conjecture to the case $k=\Q$. Assuming the geometric Lang conjecture, we reduce the weak Bombieri-Lang conjecture to $k=\mathbb{Q}$, too. 
\end{abstract}

\maketitle

In \cite{lan86}, Lang famously stated a series of conjectures centered around the scarcity of rational points of a variety of general type over fields finitely generated over $\Q$. Independently, Bombieri stated part of these conjectures. We study their reduction to number fields.

In order to describe the conjectures, let us fix some conventions and recall some definitions.

\subsection*{Conventions} Fields are tacitly assumed to be of characteristic $0$. For ease of reading, we will always use the letter $k$ for the base field, $h$ for \emph{finite} extensions of $k$ and $K$ for \emph{finitely generated} extensions of $k$. By convention, a possibly singular, non-proper variety is of general type if it is birational to a smooth, proper variety of general type.

\begin{definition*}
	A variety $X$ over $k$ is \emph{mordellic} if $X(K)$ is finite for every finitely generated extension $K$ of $k$. We stress that $X(K)$ is required to be finite for all finitely generated extensions $K$ of $k$, not only finite ones. A variety is \emph{pseudo-mordellic} if it has a non-empty open subset which is mordellic.

	A variety $X$ over $k$ is \emph{geometrically mordellic}, or \emph{GeM}, if every subvariety of $X_{\bar{k}}$ is of general type, and it is \emph{pseudo-GeM} if it has a non-empty open subset which is GeM.
\end{definition*}

Following \cite{av96}, we organize part of Lang's program in four conjectures.

\begin{description}
	\itemsep0.3em
	\item[Geometric Lang conjecture] If $X$ is a variety of general type over a field $k$ of characteristic $0$, then it is pseudo-GeM.
	\item[Lang conjecture for function fields] Let $K/k$ be a finitely generated extension of fields of characteristic $0$ with $k$ algebraically closed in $K$, and $X/K$ a variety of general type. If $X(K)$ is Zariski dense in $X$, then $X$ is birational to a variety $X_{0}$ defined over $k$ and the "non-constant points" of $X(K)\setminus X_{0}(k)$ are not Zariski-dense in $X$.
	\item[Weak Bombieri-Lang conjecture] If $X$ is a positive dimensional variety of general type over a field $k$ finitely generated over $\Q$, the set of rational points $X(k)$ is not dense.
	\item[Strong Bombieri-Lang conjecture] Every variety of general type over a field $k$ finitely generated over $\Q$ is pseudo-mordellic.
\end{description}

Faltings proved in \cite{fal94} that these conjectures hold for subvarieties of abelian varieties, but the general case remains widely open.

\begin{definition*}
	A field $k$ is \emph{weakly Lang} if, for every positive dimensional variety $X$ of general type over $k$, $X(k)$ is not dense. A field $k$ is \emph{Lang} if every positive dimensional variety $X$ of general type over $k$ is pseudo-mordellic. Clearly, Lang fields are weakly Lang, too. 
\end{definition*}

Observe that in order to check whether a field $k$ is Lang one has to consider all finitely generated extensions of $k$, while to check whether $k$ is weakly Lang it is enough to work over $k$.

Using these definitions, the weak Bombieri-Lang conjecture states that finitely generated extensions of $\Q$ are weakly Lang, while the strong version states that they are Lang.

\subsection*{Reduction to number fields}

In the last section of \cite{lan86}, Lang discusses the problem of reducing his conjectures about finitely generated extensions of $\Q$ to number fields. 

In dimension $1$, Lang cites two approaches to solve this reduction problem. The first approach is to use Mordell's conjecture for function fields, which was proved first by Manin \cite{man63} and later, with different methods, by Grauert \cite{gra65}. The second approach uses a specialization argument, see the last paragraph of \cite{lan86}. These results were available years before the Mordell conjecture for number fields was proved by Faltings in 1983 \cite{fal83}.

Let us now consider the higher dimensional case. The following easy lemmas are implicit in Lang's writings.

\begin{lemma}[Lang]\label{lfwl}
	Assume that the Lang conjecture for function fields holds, and that finite extensions of $k$ are weakly Lang. Then finitely generated extensions of $k$ are weakly Lang.
	\begin{proof}
		Follows from the definitions.
	\end{proof}
\end{lemma}

\begin{lemma}[Lang]\label{lgwl}
	Assume that the geometric Lang conjecture holds and that all finitely generated extensions of $k$ are weakly Lang. Then $k$ is Lang.
	\begin{proof}
		Let $X$ be a variety of general type over $k$. Since geometric Lang conjecture holds, there exists an open subset $U\s X$ which is GeM. We claim that $U$ is mordellic. Let $K/k$ be a finitely generated extension and define $Z\s U_{K}$ as the Zariski closure of $U_{K}(K)$. We have that $Z(K)$ is dense in $Z$ and all of its irreducible components are of general type since $U$ is GeM. Since $K$ is weakly Lang by hypothesis, then $Z$ has dimension $0$ and thus $U(K)=Z(K)$ is finite.
	\end{proof}
\end{lemma}

Combining these two lemmas, we see that the Lang conjecture for function fields (together with the geometric Lang conjecture which in dimension $1$ is trivial) represents a generalization to higher dimensions of the first approach to the reduction problem. Moreover, in the very last sentence of \cite{lan86}, Lang wishes for a generalization of the second approach (the specialization argument).

\subsection*{Our results}

Up to our knowledge, both the approaches described above are still far from being complete already in dimension $2$. We introduce a third approach based on studying the conjectures \emph{in all dimensions at the same time}. This approach reduces the strong Bombieri-Lang conjecture to $\Q$ unconditionally. Moreover, assuming the geometric Lang conjecture, we show that the strong Bombieri-Lang conjecture reduces to proving that $\Q$ is weakly Lang. This removes the need, for arithmetic purposes, of Lang's conjecture for function fields. 

\begin{TA}
	\hypertarget{TA}{Let} $K/k$ be a finitely generated extension of fields of characteristic $0$. Then $k$ is Lang if and only if $K$ is Lang. In particular, if $\Q$ is Lang then the strong Bombieri-Lang conjecture holds.
\end{TA}

\begin{TB}
	\hypertarget{TB}{Let} $k$ be a field of characteristic $0$. Assume that the geometric Lang conjecture holds. Then $k$ is Lang if and only if it is weakly Lang. 
\end{TB}

Compare \hyperlink{TB}{Theorem B} with \autoref{lgwl}: to conclude that $k$ is Lang we only need $k$ to be weakly Lang, not all of its finitely generated extensions.

Combining Theorems A and B, we immediately obtain the following.

\begin{TC}
	\hypertarget{TC}{If} the geometric Lang conjecture holds and $\Q$ is weakly Lang, then the strong Bombieri-Lang conjecture holds.
\end{TC}

We stress that both \hyperlink{TA}{Theorem A} and \hyperlink{TB}{Theorem B} take into account the various conjectures for all dimensions at the same time. For instance, if we assume the strong Bombieri-Lang conjecture \emph{only for surfaces} of general type over $\Q$, we cannot apply \hyperlink{TA}{Theorem A} to conclude that the conjecture holds for surfaces of general type over finitely generated extensions of $\Q$.

Our proofs are relatively easy, but we rely on two strong theorems: a particular case of the subadditivity of the Kodaira dimension by Viehweg \cite[Satz III]{vie82} for \hyperlink{TA}{Theorem A} and a uniformity result by Caporaso-Harris-Mazur and Abramovich-Voloch \cite{chm97}, \cite[Theorem 1.7]{av96}, \cite{abr97} for \hyperlink{TB}{Theorem B}.

Even though we are mainly interested in finitely generated extensions of $\Q$, our arguments are general. Non-standard examples are given by finitely generated extensions of $\Q(x_{1},x_{2},\dots)$: it is easy to show that if the geometric Lang conjecture and the weak Bombieri-Lang conjecture hold, then $\Q(x_{1},x_{2},\dots)$ is weakly Lang. Theorems A and B then imply that all finitely generated extensions of $\Q(x_{1},x_{2},\dots)$ are Lang.

\subsection*{Acknowledgements}

I would like to thank Dan Abramovich for suggesting some simplifications, and an anonymous referee for correcting a wrong notation.

\section{Proof of Theorem A}

\begin{lemma}\label{morext}
	Let $h/k$ be a finite extension of fields of characteristic $0$. A variety $X$ over $k$ is pseudo-mordellic if and only if $X_{h}$ is pseudo-mordellic.
	\begin{proof}
		If $U\s X$ is a mordellic open subset, then $U_{h}\s X_{h}$ is mordellic too. On the other hand, let $U\s X_{h}$ be a mordellic open subset, and let $h'$ be a Galois closure of $h/k$. Then $\bigcap_{\sigma\in\gal(h'/k)}\sigma(U_{h'})$ is Galois invariant, thus it descends to a mordellic open subset of $X$.
	\end{proof}
\end{lemma}

\begin{lemma}\label{morprod}
	Let $X_{1},\dots,X_{n}$ be varieties over a field $k$, and assume that the product $\prod_{i}X_{i}$ is pseudo-mordellic. Then each factor $X_{i}$ is pseudo-mordellic.
	\begin{proof}
		Let us prove that $X_{1}$ is pseudo-mordellic. Let $U$ be a non-empty mordellic open subset of $\prod_{i}X_{i}$. Thanks to \autoref{morext} we may assume that we have a rational point $(x_{i})_{i}\in U(k)$. Then $U\cap\left( X_{1}\times (x_{2},\dots,x_{n})\right)$ is a non-empty mordellic open subset of $X_{1}\times (x_{2},\dots, x_{n})=X_{1}$.
	\end{proof}
\end{lemma}

If $h/k$ is a finite extension and $X$ is a variety over $h$, recall that the Weil restriction $R_{h/k}(X)$ is a variety over $k$ representing the functor $S\mapsto X(S_{h})$ for schemes $S$ over $k$. The Weil restriction always exists for quasi-projective varieties, see \cite[\S 7.6]{blr90}.

If $S$ is an $h$-scheme and $\sigma:h\to h$ is a field automorphism, denote by $\sigma^{*}S$ the $h$-scheme $S\to\spec h\xar{\sigma^{\#}}\spec h$, where $\sigma^{\#}$ is the morphism of schemes induced by $\sigma$.

\begin{lemma}\label{strfin}
	Let $h/k$ be a finite extension of fields of characteristic $0$. Then $k$ is Lang if and only if $h$ is Lang.
	\begin{proof}
		If $h$ is Lang, then $k$ is Lang thanks to \autoref{morext}. Assume now that $k$ is Lang, we want to prove that $h$ is Lang. Thanks to the preceding case, we may assume that $h/k$ is Galois.
		
		Thanks to \cite[Theorem 1.3.2]{wei82}, we have that $R_{h/k}(X)_{h}\simeq \prod_{\sigma\in\gal(h/k)}\sigma^{*}X$ is a product of varieties of general type, hence $R_{h/k}(X)_{h}$ is of general type. Since we are assuming that $k$ is Lang, $R_{h/k}(X)$ is pseudo-mordellic, and thus $R_{h/k}(X)_{h}=\prod_{\sigma}\sigma^{*}X$ is pseudo-mordellic too. By \autoref{morprod}, it follows that $X$ is pseudo-mordellic.
	\end{proof}
\end{lemma}

\begin{lemma}\label{strt}
	Let $k$ be a field of characteristic $0$ and $t$ an indeterminate. If $k(t)$ is Lang, then $k$ is Lang.
	\begin{proof}
		Let $X$ be a positive dimensional variety of general type over $k$, we want to show that it is pseudo-mordellic. By hypothesis, there exists an open, mordellic non-empty subset $U\s X_{k(t)}$. Let $D\s X\times\A^{1}$ be the closure of $X_{k(t)}\setminus U$. Since $k$ is infinite and $U$ is non-empty, there exists $\lambda\in k$ such that $X\times\{\lambda\}$ is not contained in $D$ (otherwise, we would have $D=X\times\A^{1}$). Write 
		\[V=X\times\{\lambda\}\setminus D\s X\times\{\lambda\}\simeq X,\]
		we claim that $V$ is mordellic. Let $K/k$ be a finitely generated extension and $p\in V(K)$ a point, by definition of $V$ we have that $(p,\lambda)\in X\times\A^{1}$ is not contained in $D$. In particular, $\{p\}\times\A^{1}\not\s D$ and thus $p_{K(t)}\in U(K(t))$. Hence, we have that $V(K)\s U(K(t))$ is finite since $U$ is mordellic.
	\end{proof}
\end{lemma}

Let us now prove \hyperlink{TA}{Theorem A}. Let $K/k$ be a finitely generated extension. If $K$ is Lang, by \autoref{strfin} and \autoref{strt} plus an easy inductive argument on $\on{trdeg}(K/k)$ we have that $k$ is Lang.

Assume now that $k$ is Lang, we want to prove that $K$ is Lang. Thanks to \autoref{strfin} we may assume that $k$ is algebraically closed in $K$, let $V/k$ be some variety such that $k(V)=K$ and let $\eta\in V$ be its generic point. Thanks again to \autoref{strfin}, we may pass to finite extensions of $K$ and assume that $V$ is of general type.

Let $X$ be a variety of general type over $K=k(V)$, we want to prove that $X$ is pseudo-mordellic. There exists a variety $\tilde{X}$ over $k$ with a morphism $\tilde{X}\to V$ such that $\tilde{X}_{\eta}=X$. Since both $V$ and the generic fiber of $\tilde{X}\to V$ are of general type, by \cite[Satz III]{vie82} we have that $\tilde{X}$ is of general type. Since we are assuming that $k$ is Lang, there exists a non-empty mordellic open subset $U\s\tilde{X}$, i.e. $U(K')$ is finite for every finitely generated extension $K'$ of $k$. In particular, $U_{\eta}(K')\s U(K')$ is finite for every finitely generated extension $K'$ of $k(V)$ and hence $U_{\eta}$ is a non-empty mordellic open subset of $\tilde{X}_{\eta}=X$. This completes the proof of \hyperlink{TA}{Theorem A}.

\section{Proof of Theorem B}

Assume that the geometric Lang conjecture holds and that $k$ is weakly Lang, we want to show that it is Lang. Since we are assuming the geometric Lang conjecture, it is enough to show that every GeM variety over $k$ is mordellic.

Let $X/k$ be a GeM variety. We want to apply a theorem by Abramovich-Voloch \cite[Theorem 1.7]{av96}, \cite{abr97} (generalization of a result by Caporaso-Harris-Mazur \cite{chm97}) which implies that there exists an uniform bound on $|X(h)|$ for finite extensions $h/k$ of bounded degree, but to do so we need to make a small remark on their hypotheses. In fact, they assume that the base field $k$ is finitely generated over $\Q$ and that weak Bombieri-Lang conjecture holds, but a careful analysis of their proof shows that the only hypothesis on $k$ they actually use is that $k$ is weakly Lang. 

Since we are assuming that $k$ is weakly Lang, we may thus apply the uniform bound theorem by Abramovich and Voloch. Let $K/k$ be a finitely generated extension, we want to show that $X(K)$ is finite. Let $V$ be an integral scheme of finite type over $k$ whose function field is $K$, choose any generically finite, dominant rational map $V\dashrightarrow\P^n$ and let $d$ be its degree. Since rational points are dense in $\P^n$, it follows that $V_{d}=\{p\in V(\bar{k})\mid [k(p):k]\le d\}$ is dense in $V$. Thanks to the uniform bound theorem, there exists an $N$ such that $|X(h)|\le N$ for every finite extension $h/k$ with $[h:k]\le d$. Let us prove that $|X(K)|\le N$.

If by contradiction we have $N+1$ different sections $\spec K\to X$, up to shrinking $V$ we may assume that they extend to $N+1$ morphisms $f_{1},\dots,f_{N+1}:V\to X$. Since $|X(h)|\le N$ for every finite extension $h/k$ with $[h:k]\le d$, we have that for every $p\in V_{d}$ there exists a pair of different indexes $i\neq j$ with $f_{i}(p)=f_{j}(p)$. Since $V_{d}$ is dense and the pairs of different indexes are finite, there exists a pair $i\neq j$ and a subset $S\s V_{d}$ dense in $V$ such that $f_{i}(p)=f_{j}(p)$ for every $p\in S$, thus $f_{i}=f_{j}$ which gives a contradiction.

\printbibliography
	
\end{document}